\documentstyle{amsart}
\theoremstyle{plain}
\newtheorem{theorem}{Theorem}
\newtheorem{definition}[theorem]{Definition}
\newtheorem{proposition}[theorem]{Proposition}

\newtheorem{remark}[theorem]{Remark}
\newtheorem{lemma}[theorem]{Lemma}
\newtheorem{problem}[theorem]{Problem}

\newcommand{\T}{\mathbb T}
\newcommand{\C}{{\mathbb C}}
\newcommand{\D}{{\mathbb D}}

\newcommand{\cS}{{\mathcal S}}

\newcommand{\cH}{{\mathcal H}}

\numberwithin{equation}{section}
\numberwithin{theorem}{section}

\begin{document}
\begin{flushright}
\end{flushright}
\bigskip

\title[A uniqueness result]{A uniqueness result on boundary 
interpolation}

        \author{Vladimir Bolotnikov}
\address{Department of Mathematics \\
  College William and Mary \\
  Williamsburg, Virginia 23187-8795, U. S. A.}

\begin{abstract}
Let $f$ be an analytic function mapping the unit disk $\D$ to itself.
We give necessary and sufficient conditions on the local behavior of $f$
near a finite set of boundary points that requires $f$ to be a finite 
Blaschke product.
\end{abstract}

\subjclass{47A57}
\keywords{Schur functions, finite Blaschke products,
boundary interpolation problem}

\maketitle

\section{Introduction}
\setcounter{equation}{0}

The following boundary uniqueness result was presented in \cite{bk}
as an intermediate step to obtain a similar result in the multivariable 
setting of the unit ball.
\begin{theorem}
Let $f\in\cS$ and let $f(z)=z+O((z-1)^4)$ as $z\to 1$. Then $f(z)\equiv z$. 
\label{T:1.1}
\end{theorem}
Here and in what follows, $\cS$ denotes the Schur class of functions
analytic and bounded by one in modulus on the unit disk $\D$. In 
\cite{ch}, Theorem 
\ref{T:1.1} was generalized in the following way.
\begin{theorem}
Let $f\in\cS$ and let $b$ be a finite Blaschke product. Let $\tau$ be a
unimodular number and let $A_{b,\tau}=b^{-1}(\tau)=\{t_1,\ldots,t_d\}$ 
(since $b$ is a finite Blaschke product, $A_{b,\tau}$ is a finite subset 
of the unit circle $\T$). If
\begin{enumerate}
\item $f(z)=b(z)+O((z-t_1)^4)$ as $z\to t_1$ and 
\item $f(z)=b(z)+O((z-t_i)^{\ell_i})$ for some $\ell_i\ge 2$ as $z\to t_i$ 
for  $i=2,\ldots,d$, 
\end{enumerate}
then $f(z)\equiv b(z)$ on $\D$.
\label{T:1.2}
\end{theorem}
Thus, conditions in Theorem \ref{T:1.2} are sufficient to guarantee
$f(z)\equiv b(z)$. The question raised in \cite{ch} was to find 
necessary (in a sense) and sufficient conditions. The answer is given 
below. For a given real $x$, $[x]$ denotes the largest integer that does 
not exceed $x$.
\begin{theorem}
Let $f\in\cS$ and let $b$ be a finite Blaschke product of degree $d$.
Let $t_1,\ldots,t_n$ be points on $\T$ and let
\begin{equation}
f(z)=b(z)+o((z-t_i)^{m_i})\quad\mbox{for}\quad i=1,\ldots,n 
\label{1.1}
\end{equation}
as $z$ tends to $t_i$ nontangentially and where $m_1,\ldots,m_n$ are 
positive integers. If 
\begin{equation}
\left[\frac{m_1+1}{2}\right]+\ldots +\left[\frac{m_n+1}{2}\right]>d=
\deg \, b,
\label{1.2}
\end{equation}
then $f(z)\equiv b(z)$ on $\D$. Otherwise, the uniqueness result fails. 
\label{T:1.3}
\end{theorem}
In other words, the points $t_i\in\T$ can be chosen arbitrarily 
(regardless $b$) as well as degrees of convergence. To derive Theorem 
\ref{T:1.2} from  Theorem \ref{T:1.3}, note that if $\deg \, b=d$, the set 
$A_{b,\tau}$ consists of exactly $d$ points on $\T$. The assumptions in 
Theorem \ref{T:1.2} mean that \eqref{1.1} holds for $m_1=3$ and $m_i\ge 
1$ for $i=2,\ldots,d$. Then the sum on the left hand side in \eqref{1.2}
is not less than $2+(d-1)=d+1$ which is greater than $d$ and the 
result follows. The proof of Theorem \ref{T:1.3} is given in Section 4.
It relies on some recent results on boundary interpolation \cite{bk}
that we recall in Section 2 and Section 3.

\section{Boundary Schwarz-Pick matrices}
\setcounter{equation}{0}

Let $w$ be a Schur function. Then for every choice of $n\in{\mathbb N}$
and of $n$-tuples ${\bf z}=(z_1,\ldots,z_n)\in\D^n$ and ${\bf 
k}=(k_1,\ldots,k_n)\in{\mathbb
N}^n$, the {\em Schwarz-Pick matrix} $P^w_{{\bf k}}({\bf z})$
defined as
\begin{equation}
P^w_{{\bf k}}({\bf z})=\left[P^w_{k_i,k_j}(z_i,z_j)\right]_{i,j=1}^n
\label{2.1}
\end{equation}
where
\begin{equation}
P^w_{k_i,k_j}(z_i,z_j)=\left[\left.\frac{1}{\ell !r!} \,
\frac{\partial^{\ell+r}}{\partial
z^\ell\partial\bar{\zeta}^r} \, 
\frac{1-w(z)\overline{w(\zeta)}}{1-z\bar{\zeta}}
\right\vert_{{\scriptsize\begin{array}{c} z=z_i\\
\zeta=\overline{z}_j\end{array}}}
\right]_{\ell=0,\ldots,k_i-1}^{r=0,\ldots,k_j-1},
\label{2.2}
\end{equation}
is positive semidefinite. Indeed, every Schur function $w$ admits a 
de Branges--Rovnyak realization
\begin{equation}
w(z)=w(0)+zC(I_{\cH}-zA)^{-1}B\quad (z\in\D),
\label{2.3}
\end{equation}
(see \cite{dbr}) with an operator $A$ acting on an auxiliary Hilbert 
space $\cH$ and operators $B: \, \C\to \cH$ and $C: \, \cH\to\C$ such 
that the block operator ${\bf U}=\begin{bmatrix} A & B \\
C & w(0)\end{bmatrix}$ is a coisometry on $\cH\oplus \C$ (if ${\bf U}$
is unitary, representation \eqref{2.3} is called a {\em unitary 
realization} of $w$). A consequence of equality ${\bf U}{\bf U}^*=I$ is 
that
$$
\frac{1-w(z)\overline{w(\zeta)}}{1-z\bar{\zeta}}= C(I-zA)^{-1}
(I-\bar{\zeta}A^*)^{-1}C^*.
$$
Differentiating both parts in the latter identity gives
$$
\frac{1}{\ell!r!} \,
\frac{\partial^{\ell+r}}{\partial
z^\ell\partial\bar{\zeta}^r} \, 
\frac{1-w(z)\overline{w(\zeta)}}{1-z\bar{\zeta}}
=CA^\ell(I-zA)^{-\ell-1}(I-\bar{\zeta}A^*)^{-r-1}A^{*r}C^*
$$
which allows us to  represent the matrix in \eqref{2.1} as
\begin{equation}
P^w_{{\bf k}}({\bf z})=R_{{\bf k}}({\bf z})R_{{\bf k}}({\bf z})^*,
\label{2.4}
\end{equation}
where
\begin{equation}
R_{\bf k}({\bf z})=\left[\begin{array}{c}R_{k_1}(z_1)\\
\vdots \\ R_{k_n}(z_n)\end{array}\right]\quad\mbox{and}\quad
R_{k_i}(z_i)=\left[\begin{array}{c}C(I-z_iA)^{-1} \\ CA(I-z_iA)^{-2}\\
\vdots \\ CA^{k_i-1}(I-z_iA)^{-k_i}\end{array}\right],
\label{2.5}
\end{equation}  
and to conclude that $P^w_{{\bf k}}({\bf z})\ge 0$. In case when $w$
is a finite Blaschke product, the above realization arguments are more 
informative. In what follows, we will write ${\mathcal{BF}}$
for the class of all finite Blaschke products and more specifically,
${\mathcal{BF}}_d$ for the set of all Blaschke products of degree $d$.
The symbol ${\rm Dom}(w)$ will stand for the domain of holomorphy of $w$.
Apperently, the next result is well known. 
\begin{lemma}
Let $w\in{\mathcal{BF}}_d$ and let ${\bf k}=(k_1,\ldots,k_n)\in{\mathbb 
N}^n$. Then
\begin{enumerate}
\item The function $P^w_{{\bf k}}({\bf z})$ defined on $\D^n$ by formulas 
\eqref{2.1} and \eqref{2.2}, can be 
extended continuously to $({\rm Dom} (w))^n$. 
\item For every ${\bf z}\in({\rm Dom} (w))^n$, the matrix $P^w_{{\bf 
k}}({\bf z})$ is positive semidefinite and 
${\rm rank} \, P^w_{{\bf k}}({\bf z})=\min \{|{\bf k}|, \, d\}$ where   
we have set $|{\bf k}|:=k_1+\ldots+k_n$.
\end{enumerate}
\label{L:2.1}   
\end{lemma}
{\bf Proof:} Since $w$ is a rational function of McMillan degree $d$,
it admits (\cite[Chapter 4]{BGR}) a minimal realization 
\begin{equation}
w(z)=w(0)+zC(I_d-zA)^{-1}B\quad (z\in{\rm Dom} (w)),
\label{2.6}
\end{equation}
holding for all $z\in{\rm Dom} (w)$, with $\cH=\C^d$
and matrices $A\in\C^{d\times d}$, $B\in\C^{d\times 1}$,
$C\in\C^{1\times d}$ such that 
\begin{equation}
\bigcap_{j=0}^{d-1}{\rm Ker} \, CA^j=\{0\}\quad\mbox{and}\quad
\det \, (I-zA)\neq 0 \; \; (z\in{\rm Dom} (w)).
\label{2.7}   
\end{equation}
Furthermore, if $w$ inner, then the matrices $A$, $B$ and $C$ can be 
chosen so that the minimal realization \eqref{2.6} will be unitary
\cite{brkl}. The same result comes out of the de Branges--Rovnyak model:
if $w$ is inner, the de Branges--Rovnyak realization is unitary (not 
just coisometric) with the state space $\cH=H^2\ominus wH^2$; if $w\in 
{\mathcal{BF}}_d$, then $\dim\, \cH=d$ and \eqref{2.6} is obtained upon 
identifying $\cH$ with $\C^d$. 

Since realization \eqref{2.6} is unitary, formulas \eqref{2.4} and 
\eqref{2.5} hold. By \eqref{2.5}, $R_{\bf k}({\bf z})$ is analytic
on (more precisely, can be extended analytically to) $({\rm Dom} (w))^n$
and then formula \eqref{2.4} gives the desired extension of $P^w_{{\bf 
k}}({\bf z})$ to the all of $({\rm Dom}(w))^n$. By \eqref{2.5}, 
$R_{{\bf k}}({\bf z})\in\C^{|{\bf k}|\times d}$, and therefore we have 
from  \eqref{2.4} 
\begin{equation}
P^w_{{\bf k}}({\bf z})\ge 0\quad\mbox{and}\quad
{\rm rank} \, P^w_{{\bf k}}({\bf z})\le \min \{|{\bf k}|, \, d\}.
\label{2.8}
\end{equation}
On the other hand, if $|{\bf k}|=d$, the square matrix 
$R_{\bf k}({\bf z})$ is not singular. Indeed, assuming that it is 
singular, we take a nonzero vector $x\in\C^d$ such that 
$$
R_{\bf k}({\bf z})\prod_{j=1}^n (I-z_jA)^{k_j}x=0.
$$
By \eqref{2.5}, the latter matrix equation reduces to the system
of $d=|{\bf k}|$ equalities
$$
CA^{\ell}(I-z_iA)^{k_i-\ell-1}\prod_{j\neq i} (I-z_jA)^{k_j}x=0
$$
for $\ell=0,\ldots,k_i-1$ and $i=1,\ldots,n$. Expanding polynomials
leads to a homogeneous liner system (with respect to $Cx$, $CAx$,\ldots
$CA^{d-1}x$) with the nonzero Vandermonde-like determinant from which 
it follows that $GA^{\ell}x=0$ for $\ell=0,\ldots,d-1$. Then $x=0$, by the 
first relation in \eqref{2.7}, and thus, $\det \, R_{\bf k}({\bf z})\neq 
0$. By \eqref{2.4}, $P^w_{{\bf k}}({\bf z})>0$ whenever ${\bf z}\in 
({\rm Dom} (w))^n$ and $|{\bf k}|=d$. Finally if ${\bf k}=(k_1,\ldots,k_n)$ 
is any tuple with $|{\bf k}|=\widetilde{d}<d$, let $\widetilde{\bf
k}=(k_1,\ldots, k_{n-1},k_n+d-\widetilde{d})$ so that $|\widetilde{\bf 
k}|=d$. Since $P^w_{{\bf k}}({\bf z})$ is the top
$\widetilde{d}\times\widetilde{d}$ principal submatrix in  
$P^w_{\widetilde{\bf k}}({\bf z})$ and since the latter matrix is positive 
definite by the preceding analysis, we have
\begin{equation}
P^w_{{\bf k}}({\bf z})>0\quad\mbox{whenever} \; \; {\bf z}\in ({\rm 
Dom}(w))^n \; \; \mbox{and}\; \; |{\bf k}|<d.
\label{2.9}
\end{equation}
Combining \eqref{2.9} and \eqref{2.8} gives the second assertion of the
lemma and completes the proof.\qed

\medskip

Given $w\in{\mathcal{BF}}$ and a ``boundary'' tuple ${\bf 
t}=(t_1,\ldots,t_n)\in{\T}^n$, Lemma \ref{2.1} allows us to define the 
{\em boundary Schwarz-Pick matrix} $P^w_{{\bf k}}({\bf t})=
R_{{\bf k}}({\bf t})R_{{\bf k}}({\bf t})^*$ via factorization formula 
\eqref{2.4} for every ${\bf k}\in{\mathbb N}^n$. 
However, we are more interested in boundary Schwarz-Pick 
matrices for more general Schur functions. The following definition 
looks appropriate:
\begin{definition}
Given $w\in\cS$, ${\bf k}=(k_1,\ldots,k_n)\in{\mathbb N}^n$ and ${\bf 
t}=(t_1,\ldots,t_n)\in{\T}^n$, the boundary Schwarz-Pick matrix is 
defined by
\begin{equation}
P^w_{{\bf k}}({\bf t})=\lim_{{\bf z}\to{\bf t}}P^w_{{\bf k}}({\bf z})
\label{2.10}   
\end{equation}
as $z_i\in\D$ tends to $t_i$ nontangentially for $i=1,\ldots,n$, 
provided the limit in \eqref{2.10} exists. 
\label{D:2.2}
\end{definition} 
Here and in what follows, ``the limit exists'' means also that it is 
finite. By \eqref{2.1} and \eqref{2.2}, $P^w_{{\bf k}}({\bf t})$ is of the 
form
\begin{equation}
P^w_{{\bf k}}({\bf t})=\left[P^w_{k_i,k_j}(t_i,t_j)\right]_{i,j=1}^n
\label{2.11}
\end{equation}
where
\begin{equation}
P^w_{k_i,k_j}(t_i,t_j)=\lim_{{\scriptsize\begin{array}{c}z\to t_i\\ 
\zeta\to t_j\end{array}}}\left[\frac{1}{\ell !r!} \, 
\frac{\partial^{\ell+r}}
{\partial z^\ell\partial\bar{\zeta}^r} \,
\frac{1-w(z)\overline{w(\zeta)}}{1-z\bar{\zeta}}
\right]_{\ell=0,\ldots,k_i-1}^{r=0,\ldots,k_j-1}.
\label{2.12}
\end{equation}
A necessary and sufficient condition for the limits \eqref{2.12} to exist 
is that
\begin{equation}
\liminf_{z\to t_i} \frac{\partial^{2k_i-2}}{\partial
z^{k_i-1}\partial\bar{z}^{k_i-1}} \, \frac{1-|w(z)|^2}{1-|z|^2}<\infty
\quad\mbox{for}\; \; i=1,\ldots,n,
\label{2.13}
\end{equation}
where $z\in\D$ tends to $t_i$ arbitrarily (not necessarily
nontangentially). Necessity is self-evident since the bottom diagonal
entries in the diagonal blocks $P^w_{k_i,k_i}(t_i,t_i)$ are the 
nontangential (angular) limits
$$
\lim_{z,\zeta\to t_i}\frac{1}{((k_i-1)!)^2} \,
\frac{\partial^{2k_i-2}}
{\partial z^{k_i-1}\partial\bar{\zeta}^{k_i-1}} \,
\frac{1-w(z)\overline{w(\zeta)}}{1-z\bar{\zeta}} 
$$
and their existence clearly implies \eqref{2.13}. The sufficiency part 
was proved in \cite{bk} along  with some other important consequences of 
conditions \eqref{2.13} that are recalled in the following theorem.
\begin{theorem}
Let $t_1,\ldots,t_n\in\T$, $k_1,\ldots,k_n\in{\mathbb N}$, $w\in\cS$ and 
let us assume that conditions \eqref{2.13} are met. Then
\begin{enumerate}
\item The following nontangential boundary limits exist
\begin{equation}
w_j(t_i):=\lim_{z\to t_i}\frac{w^{(j)}(z)}{j!}\quad\mbox{for} \; \;
j=0,\ldots,2k_i-1; \; i=1,\ldots,n.
\label{2.14}
\end{equation}
\item The nontangential boundary limit \eqref{2.10} exists (or 
equivalently all the limits in \eqref{2.12}) exist)
and can be expressed in terms of the limits \eqref{2.14} as follows:
\begin{equation}
P^w_{k_i,k_j}(t_i,t_j)={\bf H}^w_{k_i,k_j}(t_i,t_j)
{\bf \Psi}_{k_j}(t_j){\bf T}^w_{k_j}(t_j)^*
\label{2.15}
\end{equation}
where ${\bf \Psi}_{k_j}(t_j)$ is the $k_j\times k_j$ upper triangular 
matrix with the entries
\begin{equation}
\psi_{r\ell}=\left\{\begin{array}{ccl} 0, & \mbox{if} & r>\ell \\
(-1)^\ell{\scriptsize\left(\begin{array}{c} \ell \\ r
\end{array}\right)}t_0^{\ell+r+1}, & \mbox{if} & r\leq\ell 
\end{array}\right.\quad (r,\ell=0,\ldots,k_j-1),
\label{2.16}   
\end{equation}
where ${\bf T}^w_{k_j}(t_j)$ is the lower triangular Toeplitz matrix: 
$$
{\bf T}^w_{k_j}(t_j)=\left[\begin{array}{cccc} w_0(t_j) & 0 & \ldots & 0\\
w_1(t_j) & w_0(t_j) & \ddots&\vdots \\ \vdots &\ddots&\ddots& 0\\
w_{k_j-1}(t_j) & \ldots & w_1(t_j) & w_0(t_j)\end{array}\right],
$$ 
and where ${\bf H}^w_{k_i,k_j}(t_i,t_j)$ is defined for $i=j$ as the 
Hankel matrix
\begin{equation}
{\bf H}^w_{k_i,k_i}(t_i,t_i)=\left[\begin{array}{cccc}
w_1(t_i) & w_2(t_i) & \ldots & w_{k_i}(t_j) \\
w_2(t_i) & w_3(t_i) & \ldots & w_{k_i+1}(t_i) \\
\vdots & \vdots && \vdots \\
w_{k_i}(t_i) & w_{k_i+1}(t_i) & \ldots & w_{2k_i-1}(t_i)
\end{array}\right]
\label{2.17}
\end{equation}
and entrywise (if $i\neq j$) by
\begin{eqnarray}
\left[{\bf H}(t_i,t_j)\right]_{r, \ell}
&=&\sum_{\alpha=0}^{r} (-1)^{r-\alpha}
\left(\begin{array}{c}\ell+r-\alpha \\
\ell\end{array}\right)\frac{w_{\alpha}(t_i)}
{(t_i-t_j)^{\ell+r-\alpha+1}}\nonumber \\
&&-\sum_{\beta=0}^{\ell} (-1)^{r}\left(\begin{array}{c}\ell+r-\beta \\ 
r\end{array}\right)\frac{w_{\beta}(t_j)}{(t_i-t_j)^{\ell+r-\beta+1}}
\label{2.18}
\end{eqnarray}
for $r=0,\ldots,k_i-1$ and $\ell=0,\ldots,k_j-1$.
\end{enumerate}
\begin{equation} 
(3) \; \mbox{It holds that} \; \; 
|w_0(t_i)|=1\quad(i=1,\ldots,n)\quad\mbox{and}\quad 
P^w_{{\bf k}}({\bf t})\ge 0.\qquad \qquad\quad
\label{2.19}
\end{equation}  
\label{T:2.3}  
\end{theorem}
\begin{remark}
{\rm Once the two first statements in Theorem \ref{T:2.3} are proved, the 
third statement is immediate. Inequality $P^w_{{\bf k}}({\bf t})\ge 0$ 
follows from \eqref{2.10} and the fact that $P^w_{{\bf k}}({\bf z})\ge 0$ 
for every $z\in\D$. Furthermore, existence of the limits \eqref{2.12} 
implies in particular that the nontangential boundary limits 
${\displaystyle\lim_{z\to t_i}\frac{1-|w(z)|^2}{1-|z|^2}}$ 
exist for $i=1,\ldots,n$ (and are finite) which together with existence of 
the nontangential limits $w_0(t_i)$ in \eqref{2.14} implies that 
$|w_0(t_i)|=1$}. 
\label{R:2.4} 
\end{remark} 
\begin{remark}
{\rm In case $n=1$ and $k_1=1$, Theorem \ref{T:2.3} reduces to the
classical Carath\'eodory-Julia theorem \cite{cara} on angular
derivatives.}
\label{R:2.5}
\end{remark}
\begin{remark}
{\rm In \cite{Kov}, I. Kovalishina considered the single point case ($n=1$
and $k_1>1$) under an additional assumption that $w$ satisfies
the symmetry relation $w(z)\overline{w(1/\bar{z})}\equiv 1$ in some
neighborhood of $t_1$. A  remarkable ``Hankel-${\bf \Psi}$-Toeplitz''
structure \eqref{2.15} of $P^w_{k_1,k_1}(t_1,t_1)$ has been observed
there}.
\label{R:2.6}
\end{remark}
Carath\'eodory-Julia type conditions \eqref{2.13} are worth a 
formal definition.
\begin{definition}
Given $n$-tuples ${\bf t}=(t_1,\ldots,t_n)\in{\T}^n$ and
${\bf k}=(k_1,\ldots,k_n)\in{\mathbb N}^n$, a Schur function $w$ is 
said to belong to the class $\cS_{\bf k}({\bf t})$ if it meets
conditions \eqref{2.13}.
\label{D:2.7}
\end{definition}
Statement (1) in Theorem
\ref{T:2.3} shows that the definition \eqref{2.10} of the boundary
Schwarz-Pick matrix $P^w_{{\bf k}}({\bf t})$ makes sense if and only if 
$w\in\cS_{\bf k}({\bf t})$. Statement (2) expresses $P^w_{{\bf 
k}}({\bf t})$ in terms of boundary limits of $w$ and of its
derivatives. An interesting point in \eqref{2.19} is that 
the block matrix $P^w_{{\bf k}}({\bf t})$ of the form \eqref{2.11} 
constructed  via structured formulas \eqref{2.15}--\eqref{2.18} (rather 
than by the limits \eqref{2.12}) does not look like a Hermitian matrix 
and nevertheless, it turns out to be Hermitian (and even positive 
semidefinite) due to conditions \eqref{2.13}. The next theorem (see 
\cite{bk} for the proof) shows that relations \eqref{2.19} are 
characteristic for the class $\cS_{\bf k}({\bf t})$.
\begin{theorem}
Let $w$ be a Schur function, let ${\bf t}\in{\T}^n$, ${\bf k}\in{\mathbb 
N}^n$ and let us assume that the nontangential limits \eqref{2.14} 
exist and satisfy conditions \eqref{2.19} where $P^w_{{\bf k}}({\bf t})$
is the matrix constructed from the limits \eqref{2.14} via formulas 
\eqref{2.15}--\eqref{2.18}. Then $w\in\cS_{\bf k}({\bf t})$. 
\label{T:2.7}
\end{theorem}
From the computational point of view, it is much easier to construct
the boundary Schwarz-Pick matrix $P^w_{{\bf k}}({\bf t})$ via formulas
\eqref{2.15}--\eqref{2.18}, than by \eqref{2.12} (for example, if $w$ is 
a rational function, the boundary limits $w_i(t_j)$ are just the Taylor 
coefficients from the expansion of $w$ around $t_i$). However, as follows 
from Theorems \ref{T:2.3} and \ref{T:2.7}, the matrix constructed in this 
way will be indeed the boundary Schwarz-Pick matrix if and only if 
conditions \eqref{2.19} are satisfied.

\section{Boundary interpolation for classes $\cS_{\bf k}({\bf t})$}
\setcounter{equation}{0}

The following interpolation problem has been studied in \cite{bk}.
\begin{problem}
Given ${\bf t}=(t_1,\ldots,t_n)\in{\T}^n$,
${\bf k}=(k_1,\ldots,k_n)\in{\mathbb N}^n$ and numbers $b_{ij}\in\C$ 
$(j=0,\ldots,k_i-1; \; i=1,\ldots,n)$, find all functions $f\in\cS_{\bf 
k}({\bf t})$ such that
\begin{equation}
f_j(t_i):=\lim_{z\to t_i}\frac{f^{(j)}(z)}{j!}=b_{ij}\quad
(j=0,\ldots,2k_i-1; \; i=1,\ldots,n).
\label{3.1}  
\end{equation}
where all the limits are nontangential.
\label{P:3.1}
\end{problem}
This interpolation problem makes perfect sense: if $f$ belongs 
to $\cS_{\bf k}({\bf t})$, then the nontangential limits in \eqref{3.1}
exist by Theorem \ref{T:2.3}; we just want them to be equal to 
the preassigned numbers. Let define the $|{\bf k}|\times |{\bf k}|$ 
matrix $P$ (the {\em Pick matrix} of the problem) by formulas similar
to \eqref{2.15}--\eqref{2.18}, but with $w_j(t_i)$ replaced by $b_{ij}$:
\begin{equation}
P=\left[P_{ij}\right]_{i,j=1}^n\quad\mbox{with}\quad
P_{ij} =H_{ij}\cdot {\bf \Psi}_{k_j}(t_j)\cdot T_j^*,
\label{3.3}
\end{equation}
where ${\bf \Psi}_{k_j}(t_j)$ is the upper triangular matrix with the 
entries given in \eqref{2.16}, where $T_i$ is the lower triangular 
Toeplitz matrix and $H_{ii}$ is the Hankel matrix defined by
\begin{equation}
T_i=\left[\begin{array}{ccc} b_{i,0} & & 0 \\ 
\vdots & \ddots & \\ 
b_{i,k_j-1} & \ldots & b_{i,0}\end{array}\right],\quad
H_{ii}=\left[\begin{array}{ccc} b_{i,1} & \cdots &
b_{i,k_i}\\ \vdots & & \vdots\\ b_{i,k_i} & \cdots &
b_{i,2k_i-1}\end{array}\right]
\label{3.4}
\end{equation}
for $i=1,\ldots,n$ and where the matrices $H_{ij}$ (for $i\neq j$)
are defined entrywise by
\begin{eqnarray}
\left[H_{ij}\right]_{r, \ell}
&=&\sum_{\alpha=0}^{r} (-1)^{r-\alpha}  
\left(\begin{array}{c}\ell+r-\alpha \\
\ell\end{array}\right)\frac{b_{i,\alpha}}
{(t_i-t_j)^{\ell+r-\alpha+1}}\nonumber \\
&&-\sum_{\beta=0}^{\ell} (-1)^{r}\left(\begin{array}{c}\ell+r-\beta \\
r\end{array}\right)\frac{b_{j,\beta}}{(t_i-t_j)^{\ell+r-\beta+1}}.
\label{3.2}
\end{eqnarray}
The purpose of this construction is clear: the matrix $P$ constructed 
above depends on the interpolation data only; on the other hand, for every
solution $f$ of Problem \ref{P:3.1}, the boundary Schwarz-Pick matrix
$P^f_{{\bf k}}({\bf t})$ must be equal to $P$, by the very construction.
\begin{theorem}
Let $P$ be the matrix defined in \eqref{3.3}. Then
\begin{enumerate}
\item If Problem \ref{P:3.1} has a solution, then 
\begin{equation}
|b_{i,0}|=1\quad(i=1,\ldots,n)\quad\mbox{and}\quad P\ge 0.
\label{3.5}   
\end{equation}
\item If \eqref{3.5} holds and $P>0$, then Problem \ref{P:3.1} 
has  infinitely many solutions.
\item If \eqref{3.5} holds and $P$ is singular, then Problem 
\ref{P:3.1} has at most one solution.
\item If \eqref{3.5} holds and $f$ is a Schur function satisfying 
conditions \eqref{3.1}, then necessarily $f\in\cS_{\bf k}({\bf t})$.
\end{enumerate}
\label{T:3.2}
\end{theorem}
The first statement follows from Statement (3) in Theorem \ref{T:2.3},
since $b_{i,0}=f_0(t_i)$ and $P^f_{{\bf k}}({\bf t})=P$ for every solution 
$f$ of Problem \ref{P:3.1}. The last statement follows from Theorem
\ref{T:2.7}. The second statement is proved in \cite{bk}
where moreover, a linear fractional parametrization of all solutions 
of Problem \ref{P:3.1} (in case $P>0$) is given. The third statement 
also was proved in \cite{bk}. 

\smallskip

The proof of Theorem \ref{T:1.3} will rest on Theorem \ref{T:3.2} and 
the following simple observation.
\begin{proposition}
Let $\widetilde{P}=[p_{ij}]$ be an $r\times r$ Hermitian matrix and let 
us assume that its principal 
submatrix $P=[p_{i_\alpha,i_\beta}]_{\alpha,\beta=1}^{\ell}$
is positive definite. Then $\widetilde{P}$ can be turned into a positive 
definite  matrix upon an appropriate modification of the 
$r-\ell$ diagonal entries $p_{ii}$ for $i\not\in\{i_1,\ldots,i_{\ell}\}$
(we will call these entries the diagonal entries of $\widetilde{P}$
complementary to the principal submatrix $P$).
\label{P:3.3}
\end{proposition}
{\bf Proof:} Without loss of generality we can assume that $P$ is the 
leading principal submatrix of $\widetilde{P}$ so that  
$\widetilde{P}=\left[\begin{array}{cc} P & 
R^* \\ R & D\end{array}\right]$. Let us modify the diagonal entries in 
$D$ as follows: 
$$
\widetilde{P}^\prime=\left[\begin{array}{cc} P &
R^* \\ R & D^\prime\end{array}\right] \quad \mbox{where} \; \; D^\prime=
D+\rho I_{r-\ell} \; \; (\rho >0).
$$
Since $P>0$, the factorization formula
$$
\left[\begin{array}{cc} P & R^* \\ R & D^\prime\end{array}\right]=   
\left[\begin{array}{cc} I_\ell & 0 \\ RP^{-1} & 
I_{r-\ell}\end{array}\right]
\left[\begin{array}{cc} P & 0 \\ 0 & D^\prime-RP^{-1}R^*\end{array}\right]
\left[\begin{array}{cc} I_\ell & P^{-1}R^* \\ 0 & 
I_{r-\ell}\end{array}\right]
$$
shows that $\widetilde{P}^\prime>0$ if and only if 
$D^\prime-RP^{-1}R^* =\rho I_{r-\ell}+D-RP^{-1}R^*>0$ and the latter 
inequality indeed can be  achieved if $\rho$ is large enough.\qed

\section{Proof of Theorem \ref{T:1.3}}
\setcounter{equation}{0}

Let us assume for a moment that the Schur function $f$ in \eqref{1.1} is 
not given and let us consider the following interpolation problem.
\begin{problem} 
Given $t_1,\ldots,t_n\in\T$, $m_1,\ldots,m_n\in{\mathbb N}$ and 
$b\in{\mathcal{BF}}_d$, find all Schur functions
$f$ satisfying asymptotic equations \eqref{1.1}.
\label{P:4.1}
\end{problem}
Note that conditions \eqref{1.1} can be reformulated equivalently
(see e.g., \cite[Corollary 7.9]{boldym1} for the proof) as follows:
{\em the nontangential boundary limits $f_j(t_i)$ exist and satisfy
\begin{equation}
f_j(t_i):=\lim_{z\to 
t_i}\frac{f^{(j)}(z)}{j!}=\frac{b^{(j)}(t_i)}{j!}=:b_{ij}\quad \mbox{for} 
\; \; j=0,\ldots,m_i; \; i=1,\ldots,n.
\label{3.6}
\end{equation}}
Define the integers $k_i:=\left[\frac{m_i+1}{2}\right]$ for $i=1,\ldots,n$
so that $m_i=2k_i-1$ or $m_i=2k_i$. Reindexing if necessary, we can assume 
without loss of generality that the first $\ell$ integers 
$m_1,\ldots,m_\ell$ are odd while the remaining ones (if any) are even.
Now we split conditions \eqref{3.6} into two parts:
\begin{equation}
f_j(t_i)=\frac{b^{(j)}(t_i)}{j!}=:b_{ij}\quad \mbox{for}
\; \; j=0,\ldots,2k_i-1; \; i=1,\ldots,n
\label{3.7}   
\end{equation}
and 
\begin{equation}
f_{2k_i}(t_i)=\frac{b^{(2k_i)}(t_i)}{(2k_i)!}=:b_{i,2k_i}\quad 
\mbox{for} \; \; i=\ell+1,\ldots,n.
\label{3.8}
\end{equation} 
First we consider the interpolation problem with interpolation conditions
\eqref{3.7} (this problem is ``truncated'' with respect to
Problem \ref{P:4.1}). This problem looks like Problem \ref{P:3.1};
however, it is more special, since that data $\{b_{ij}\}$ comes from 
certain $b\in{\mathcal {BF}}_d$. In other words, the Pick 
matrix $P$ of the problem \eqref{3.7} coincides with the boundary 
Schwarz-Pick matrix $P^b_{\bf k}({\bf t})$. Then we may conclude by Lemma 
\ref{L:2.1} that $P\ge 0$ and 
\begin{equation}
{\rm rank} \, P=\min \{|{\bf k}|, \, d\}.
\label{3.9}
\end{equation}
Thus, the second condition in \eqref{3.5} is met, while the first 
condition  holds since $b_{i,0}=b(t_i)$ and $b\in{\mathcal {BF}}$. 
Assuming that inequality \eqref{1.2} is in force, i.e., that 
$$
d<\sum_{i=1}^n\left[\frac{m_i+1}{2}\right]=\sum_{i=1}^nk_i=|{\bf k}|
$$
we conclude from \eqref{3.9} that $P$ is singular and then by Statement 
(3) in Theorem \ref{T:3.2}, there is at most one $f\in\cS$ satisfying
conditions \eqref{3.7}. Therefore (since \eqref{3.7} is just part of 
\eqref{3.6}), there is at most one $f\in\cS$ satisfying
conditions \eqref{3.6}. A self-evident observation that the Schur function
$b$ does satisfy \eqref{3.6} (this information is contained in 
\eqref{3.6}) gives the desired uniqueness: there are no functions $f$
in $\cS$ different from $b$ that satisfy interpolation 
conditions \eqref{3.6} or, equivalently, asymptotic equalities 
\eqref{1.1}. Thus, once \eqref{1.2} is in force and $f$ is subject to
\eqref{1.1}, we have necessarily $f(z)\equiv b(z)$. This completes the 
proof of the first statement in Theorem \ref{T:1.3}. It remains to show 
that the uniqueness result fails whenever $|{\bf k}|\le d$. In this case 
we conclude from \eqref{3.9} that the $|{\bf k}|\times |{\bf k}|$ matrix 
$P$ is positive definite and then, by Statement (2) in Theorem 
\ref{T:3.2}, there are infinitely many Schur functions $f$ satisfying
conditions \eqref{3.7}. In case all $m_i$'s are odd, this completes the 
proof: conditions  \eqref{3.7} are identical with \eqref{3.6} and thus, 
there  are infinitely many Schur functions satisfying asymptotic \eqref{1.1}.
The general case (when the set of conditions \eqref{3.8} is not empty)
requires one step more.

\smallskip

Assuming that $|{\bf k}|\le d$ so that the Pick matrix 
$P=P^b_{\bf k}({\bf t})$ corresponding to interpolation problem 
\eqref{3.7} is positive definite and that $\ell<n$ in 
\eqref{3.8}, let us attach interpolation 
conditions 
\begin{equation}
f_{2k_i+1}(t_i)=\frac{b^{(2k_i+1)}(t_i)}{(2k_i+1)!}=:b_{i,2k_i+1}\quad
\mbox{for} \; \; i=\ell+1,\ldots,n
\label{3.10}
\end{equation}
to \eqref{3.8} and let us consider the extended interpolation problem (for 
Schur class functions) with interpolation conditions \eqref{3.7}, 
\eqref{3.8} and \eqref{3.10}. The collection of $b_{ij}$'s appearing in   
\eqref{3.7} and \eqref{3.8} will be called the {\em original data}, 
the collection $\{b_{i,2k_i+1}\}$ from \eqref{3.10} will be called the 
{\em supplementary data} whereas their union will be referred to as to the 
{\em extended data}.

\smallskip

For the extended interpolation problem we have an even number of 
conditions for each interpolating point $t_i$ which allows us to 
construct the corresponding extended Pick matrix $\widetilde{P}$
via formulas \eqref{3.3}:
\begin{equation}
\widetilde{P}=\left[\widetilde{P}_{ij}\right]_{i,j=1}^n
\quad\mbox{where}\quad
\widetilde{P}_{ij} =\widetilde{H}_{ij}\cdot {\bf
\Psi}_{\widetilde{k}_j}(t_j)\cdot \widetilde{T}_j^*
\label{3.11}
\end{equation}
and where $\widetilde{H}_{ij}$ and $\widetilde{T}_j$ are defined by 
formulas \eqref{3.4}, \eqref{3.2} with $k_i$ replaced by 
$\widetilde{k}_i$. It is clear that $\widetilde{P}$ coincides with
the boundary Schwarz-Pick  matrix $P^b_{\widetilde{\bf k}}({\bf t})$
based on the same $b\in{\mathcal {BF}}_d$, the same ${\bf
t}=(t_1,\ldots,t_n)\in{\T}^n$ and 
$$
\widetilde{\bf k}=(\widetilde{k}_1,\ldots,\widetilde{k}_n)=
(k_1,\ldots,k_\ell,k_{\ell+1}+1,\ldots,k_n+1)\in{\mathbb N}^n.
$$
Of course, all the entries in $\widetilde{P}$ are expressed in terms of
the extended data. However, it turns out that all its entries but
$\ell$ diagonal ones are uniquely determined from the {\em original data}.   
Indeed, if $i\neq j$, then $\widetilde{H}_{ij}$ and $\widetilde{T}_j$
(and therefore, $\widetilde{P}_{ij}$) are expressed via formulas 
\eqref{3.4}, \eqref{3.2} in terms  the numbers 
$b_{i,0},\ldots,b_{i,\widetilde{k}_i-1}$ and 
$b_{j,0},\ldots,b_{j,\widetilde{k}_j-1}$ all of which are contained in 
the original data, since $\widetilde{k}_i-1\le k_i\le 2k_i-1$.

\smallskip

Now we examine the diagonal blocks $\widetilde{P}_{ii}$ for $i>\ell$
(if $i\le \ell$, then $\widetilde{P}_{ii}=P_{ii}$ is completely determined
by the original data). By \eqref{3.11} and \eqref{3.4},
\begin{equation}
\widetilde{P}_{ii}=\left[\begin{array}{cccc} b_{i,1} & \cdots & b_{i,k_i}&
b_{i,k_i+1}\\ \vdots & & \vdots &\vdots\\
b_{i,k_i} & \ldots & b_{i,2k_i-1}& b_{i,2k_i}\\
b_{i,k_i+1} & \ldots & b_{i,2k_i} & b_{i,2k_i+1}\end{array}\right]
{\bf \Psi}_{k_i+1}(t_i)\left[\begin{array}{ccc}\overline{b}_{i,0} &
\ldots & \overline{b}_{i,k_i}\\ & \ddots & \vdots \\
0 & &  \overline{b}_{i,0} \end{array}\right].
\label{3.12}
\end{equation}
It is readily seen from \eqref{3.12} that the only entry in
$\widetilde{P}_{ii}$ that depends on the supplementary data
is the the bottom diagonal entry
\begin{equation}
\gamma_i:=\left[\widetilde{P}_{ii}\right]_{k_i,k_i}
=\left[\begin{array}{ccc}b_{i,k_i+1} & \cdots &
b_{i,2k_i+1}\end{array}\right]{\bf \Psi}_{k_i+1}(t_i)
\left[\begin{array}{ccc} b_{i,k_i}&\cdots &
b_{i,0}\end{array}\right]^*
\label{3.13}
\end{equation}
which, on account of \eqref{2.16}, can be written as
\begin{eqnarray}
\gamma_i&=&(-1)^{k_i}t_i^{2k_i+1}b_{i,2k_i+1}\overline{b}_{i,0}\label{3.14}\\
&&+\sum_{r=0}^{k_i-1}b_{i,k_i+r+1}\sum_{j=k_i+r}^{k_i}
(-1)^{j}t_i^{k_i+r+j+1}\left(\begin{array}{c} j \\ k_i+r
\end{array}\right)\overline{b}_{i,n_i-j}.\nonumber
\end{eqnarray}
Since $\widetilde{P}$ coincides with
the boundary Schwarz-Pick  matrix $P^b_{\widetilde{\bf k}}({\bf t})$,
it is positive semidefinite (by Lemma \ref{L:2.1}) and Hermitian, in 
particular. Furthermore, the Pick matrix $P=P^b_{{\bf  
k}}({\bf t})$ of the problem \eqref{3.7} is a positive definite principal
submatrix of $\widetilde{P}$. The diagonal entries in $\widetilde{P}$
complementary to $P$ are exactly $\gamma_i$'s from \eqref{3.13}, the
bottom diagonal entries in the blocks $\widetilde{P}_{ii}$ of
$\widetilde{P}$ for $i=\ell+1,\ldots,n$. By Proposition \ref{P:3.3}, upon
replacing $\gamma_i$ in $\widetilde{P}$ by appropriately chosen
(sufficiently large) positive numbers $\gamma_i^\prime$ (for
$i=\ell+1,\ldots,n$) and keeping all the other entries the same, one gets
a positive definite matrix $\widetilde{P}^\prime$. Furthermore, for each 
chosen $\gamma_i^\prime$, there exists (the unique) $b_{i,2k_i+1}^\prime$ 
such that
\begin{eqnarray}
\gamma_i^\prime&=&(-1)^{k_i}t_i^{2k_i+1}b_{i,2k_i+1}^\prime\overline{b}_{i,0}
\nonumber\\
&&+\sum_{r=0}^{k_i-1}b_{i,k_i+r+1}\sum_{j=k_i+r}^{k_i}
(-1)^{j}t_i^{k_i+r+j+1}\left(\begin{array}{c} j \\ k_i+r
\end{array}\right)\overline{b}_{i,n_i-j}\nonumber
\end{eqnarray}
(since  ${b}_{i,0}\neq 0$, the latter equality can be solved for
$b^\prime_{i,2k_i+1}$). Now we replace the supplementary interpolation 
conditions \eqref{3.10} by
\begin{equation}
f_{2k_i+1}(t_i)=b^\prime_{i,2k_i+1}\quad
\mbox{for} \; \; i=\ell+1,\ldots,n
\label{3.15}
\end{equation}
where the numbers on the right have nothing to do with the finite
Blaschke product $b$ anymore. It is easily seen that the Pick matrix of 
the modified extended interpolation problem with interpolation 
conditions \eqref{3.7}, \eqref{3.8} and \eqref{3.15} is 
$\widetilde{P}^\prime$. Since it is positive definite, there are (by 
Statement (2) in Theorem \ref{T:3.2}) infinitely many Schur functions $f$ 
satisfying these interpolation conditions. Thus, there are infinitely many 
Schur functions satisfying \eqref{3.7}, \eqref{3.8} (that is, \eqref{3.6}) 
or equivalently, the asymptotic equalities \eqref{1.1}. Thus, the uniqueness
conclusion in Theorem \ref{T:1.3} fails which completes the proof.\qed


\end{document}